\documentclass[12pt,a4paper]{article}
\usepackage{amsfonts,amsmath,amssymb,indentfirst}

\newcommand{\be}{\mathbb{B}}
\newcommand{\er}{\mathbb{R}}

\newcommand{\e}{\epsilon}

\newtheorem{Teorema}{Theorem}[section]
\newtheorem{Propriedade}[Teorema]{Proposition}
\newtheorem{Rem}[Teorema]{Remark}

\newtheorem{Lema}[Teorema]{Lemma}

\begin{document}
\title{On a quasilinear non-local Benney System}
\author{Jo\~ao-Paulo Dias and Filipe Oliveira}
\date{}
\maketitle
\begin{abstract}
 We study the quasilinear non-local Benney System 
 \begin{displaymath}
 \left\{\begin{array}{llll}
iu_t+u_{xx}=|u|^2u+buv
\\
\displaystyle v_t+a\Big(\int_{\er^+}v^2dx\Big)v_x=-b(|u|^2)_x,\quad (x,t)\in\er^+\times [0,T],\, T>0.
        \end{array}\right.
\end{displaymath}
We establish the existence and uniqueness of strong local solutions to the corresponding Cauchy problem and show, under certain conditions, the blow-up of such solutions in finite time. Furthermore, we prove the existence of global weak solutions and exhibit bound-state solutions to this system.

\end{abstract}
\section{Introduction}
The Benney systems (\cite{intro1}) appear as universal models for describing the interaction between unidimensional short and long waves propagating in dispersive media. For instance, in water waves general theory, such systems can be used to represent the interaction of capillary and gravity waves or internal and surface waves (\cite{intro3,intro8}). In plasma physics, the Benney systems are known to accurately model the interaction between Alfv\'en and Langmuir waves. Other applications include the study of resonant waves in geometric optics (\cite{intro9}), the long wave/short wave interactions in bubbly liquids (\cite{intro1bis}) or optical-microwave interactions in nonlinear mediums (\cite{intro14bis}).\\

In the case where the amplitudes of the short and long waves are comparable, long waves become considerably weaker and the following quasilinear version of the Benney systems should then be considered (see for instance Eqs. (3.27) and (3.28) in \cite{intro1}):
\begin{equation}
 \label{quasi}
\left\{\begin{array}{lll}
 iu_t+u_{xx}=|u|^2u+buv\\
 v_t+(f(v))_x=c(|u|^2)_x.
 \end{array}\right.
\end{equation}
Here, $b$ and $c$ are real constants, $u=u_{(y)}+iu_{(z)}$ denotes, in complex notation, the transverse components of the short wave with respect to the direction of propagation $(0x)$, $v$ is the perturbation induced by the long wave and $f$ is a nonlinear polynomial. Contrarely to the case where $f$ is linear (which corresponds to the classical Benney systems), only a few mathematical results concerning \eqref{quasi} exist in the litterature (see \cite{DF} and \cite{DFF} for the existence of weak solutions, \cite{DFO} and \cite{intro2bis} for the study of the Cauchy problem in $\er$ and in the half-line, and \cite{intro4bis} for the study of shock wave solutions).\\

Recently, in order to study nonlinear ion transport phenomena, some attention has been devoted to non-local versions of \eqref{quasi} (see \cite{CPR}, \cite{MH}). This constitutes the motivation to consider the following nonlocal quasilinear evolution system, which can be used to model the transport of ions in an electrolytic solution in a large tube:
\begin{equation}
 \label{sistemain}
 \left\{\begin{array}{llll}
iu_t+u_{xx}=|u|^2u+buv&(i)
\\
\displaystyle v_t+a\Big(\int_{\er^+}v^2dx\Big)v_x=-b(|u|^2)_x.&(ii)
        \end{array}\right.
\end{equation}
This system describes the interaction between electrons and ions, whose density and velocity are represented respectively by
$u$ (complex-valued) and $v$ (real-valued). Here, $(x,t)\in [0,+\infty[\times [0,T]$, $a>0$ and $b\in\er$. We will consider initial data
\begin{equation}
\label{idata}
 u(x,0)=u_0(x)\in H_0^1(\er^+)\cap H^3(\er^+),\,v(x,0)=v_0(x)\in H_0^1(\er^+)\cap H^2(\er^+)
\end{equation}
and boundary conditions
\begin{equation}
 \label{bc}
 \forall t,\,u(0,t)=v(0,t)=\lim_{x\to +\infty} u(x,t)=\lim_{x\to +\infty} v(x,t)=0.
\end{equation}

The main difficulty of \eqref{sistemain} lies  in the presence of the nonlocal term\\ $a\Big(\int_{\er^+}v^2dx\Big)v_x$. To the best of our knowledge, results concerning this system do not seem available in the literature. It is easy to check that if $b=0$ the Cauchy problem of the (decoupled) system can be easily solved, since the transport equation becomes linear ($\int_{\er_+}v^2(x,t)dx=\int_{\er_+}v_0(x)dx$). Hence, in the present paper, we will consider the case $b\neq 0$.

\bigskip

The rest of this paper is organized as follows:\\

In the next Section we prove the existence and uniqueness of a local in time strong solution to the Cauchy Problem (\ref{sistemain}-\ref{idata}-\ref{bc}). This will be achieved by using a technique based on a T. Kato's result (see \cite{K}) concerning hyperbolic evolution systems, which we will adapt to the non-local case.The derivative-loss in the right-hand-side of (\ref{sistemain}-i) will be handled by performing an adequate change of the dependent variables, following the ideas in \cite{DFO} and \cite{O}.

\bigskip

In Section 3, we will start by exhibiting a few invariants for the flow of \eqref{sistemain} and, with some additional requirements on the initial data, we derive a blow-up in time result of the local strong solution of the Cauchy problem  (\ref{sistemain}-\ref{idata}-\ref{bc}). We will apply a virial technique developped in the seminal works of R.T. Glassey (cf. \cite{G}) concerning nonlinear Schr\"odinger equations (see also \cite{intro2bis} for a related result in the local situation). The function
$$t\to \int |x|^2|u(x,t)|^2dt,$$
used in the classical situation, will be replaced by
$$t\to \frac 12\int x^2|u|^2dx+\int_0^t\int xv^2dxd\tau+b^2\int_0^t\int x|u|^2dxd\tau,$$
for which a suitable inequality will be proved.

\bigskip

In Section 4 we prove  the existence of a global in time (that is, in every time interval of the form $[0,T]$, $T>0$) weak solution of the Cauchy problem (\ref{sistemain}-\ref{idata}-\ref{bc}) by applying the vanishing viscosity method to system (\ref{sistemain}) with a parabolic regularization of the transport equation and by deriving suitable {\it a priori} estimates leading to a term of boundary-layer type.

\bigskip

Finally, in the last Section we show, in some special cases, the existence of bound-state solutions to \eqref{sistemain}.

\section{Local existence and uniqueness}
Following \cite{DFO}, for $u_0\in H_0^1(\er^+)\cap H^3(\er^+)$ and $v_0\in H_0^1(\er^+)\cap H^2(\er^+)$ let us take $(u,v)$ a possible solution in $\er^+\times [0,T]$, $T>0$, of the Cauchy problem (\ref{sistemain}-\ref{idata}-\ref{bc}) and make the following formal computations:\\

\noindent
By setting $F=u_t$ we obtain, from \eqref{sistemain},
$iF+u_{xx}-u=|u|^2u+buv-u$
and so, with $\displaystyle\Delta=\frac{\partial^2}{\partial x^2}$, 
\begin{equation}
	\label{alt1}
	u=(\Delta-1)^{-1}(|u|^2u+u(bv-1)-iF).
\end{equation}
Differentiating (\ref{sistemain}-i) with respect to $t$ leads to
	$$iF_t+F_{xx}=2|u|^2F+u^2\overline{F}+bFv+buv_t$$
and, using (\ref{sistemain}-ii),
\begin{equation}
\label{alt2}
iF_t+F_{xx}=2|u|^2F+u^2\overline{F}+bFv-b^2u(|u|^2)_x-abu\Big(\int v^2dx\Big)v_x,
\end{equation}
where $\displaystyle\int v^2dx=\int_{\er^+}v^2(x,t)dx.$

\bigskip

These formal computations suggest us to consider the following Cauchy problem:
\begin{equation}
\label{sistaltf}
\left\{\begin{array}{llll}
iF_t+F_{xx}=2|u|^2F+u^2\overline{F}+bFv-b^2u|\tilde{u}|^2_x-abu\Big(\int v^2dx\Big)v_x&(i)\\
v_t+a(\int v^2dx\Big)v_x=-b|\tilde{u}|^2_x&(ii),\\
\end{array}\right.
\end{equation}
\begin{equation}
\label{idataf}
F(x,0)=F_0(x)\in H_0^1(\er^+), \quad v(x,0)=v_0(x)\in H^2(\er^+)\cap H_0^1(\er^+),
\end{equation}
where $u$ and $\tilde{u}$ are given in terms of $F$ by
$$u(x,t)=u_0(x)+\int_0^t F(x,\tau)d\tau$$
and
$$\tilde{u}(x,t)=(\Delta-1)^{-1}(|u|^2u+u(bv-1)-iF).$$
Since, for $F\in H_0^1$,  the regularization provided by the operator $(\Delta-1)^{-1}$ puts $\tilde{u}$ in $H^1_0\cap H^3$ and there are no derivative losses on the right-hand-side of (\ref{sistaltf}-i).\\

The following lemma will be proved using an adaptation of the general result of T. Kato (Theorem 6 in \cite{K}) on quasilinear systems:
\begin{Lema}
\label{lema1}
Let $(F_0,v_0)\in H^1_0(\er^+)\times(H^2(\er^+)\cap H^1_0(\er^+))$.\\
Then, there exists $T>0$ and a unique strong solution $(F,v)$ of the Cauchy problem (\ref{sistaltf}-\ref{idataf}), with
\begin{multline}
(F,v)\in \Big(C^j([0,T];H^{1-2j}(\er^+))\times C^j([0,T];H^{2-j}(\er^+))\Big)\bigcap\\
\Big(C([0,T];H^1_0(\er^+))\times C([0,T];H^{1}_0(\er^+))\Big),\quad j=0,1.
\end{multline}
\end{Lema}
Note that this Lemma implies the following result:
\begin{Teorema}
\label{TEL}
Assume that $(u_0,v_0)$ verifies \eqref{idata}. Then there exists $T>0$ and a unique strong solution $(u,v)$ of the Cauchy problem (\ref{sistemain}-\ref{idata}-\ref{bc}) with
\begin{multline}
(u,v)\in \Big(C^j([0,T];H^{3-2j}(\er^+))\times C^j([0,T];H^{2-j}(\er^+))\Big)\bigcap\\
\Big(C^j([0,T];H^1_0(\er^+))\times C^j([0,T];H^{1}_0(\er^+))\Big),\quad j=0,1.
\end{multline} 
\end{Teorema}

Indeed, if $(F,v)$ is a solution of (\ref{sistaltf}), we obtain $u_t=F$ and $u(x,0)=u_0(x)$. We derive
\begin{multline*}
(iu_t+u_{xx})_t=2|u|^2F+u^2\overline{F}+bFv-b^2u|\tilde{u}|^2_x-abu\Big(\int v^2dx\Big)v_x\\
=2|u|^2u_t+u^2\overline{u}_t+bu_tv+buv_t.
\end{multline*}
Hence, $(iu_t+u_{xx}-|u|^2u-buv)_t=0$ and $iu_t+u_{xx}-|u|^2u-buv=\phi_0(x)$, where
$$\phi_0=iF_0+u_0''-|u_0|^2u_0-bu_0v_0.$$
If we set $F_0=i(u_0''-|u_0|^2u_0-bu_0v_0)$, we obtain $\phi_0=0$ and $(u,v)$ satisfies (\ref{sistemain}-i).\\
Moreover,
\begin{equation}
\label{26}
u=(\Delta-1)^{-1}(|u|^2u+u(bv-1)-iu_t)
\end{equation}
and therefore $\tilde{u}=u$ and $(u,v)$ satisfies (\ref{sistemain}-ii).\\
Finally, note that $u_t=F\in C([0,T],H_0^1(\er^+))$ and that
$$u(\cdot,t)=u_0(\cdot)+\int_0^t F(\cdot,\tau)d\tau\in C([0,T];H_0^1(\er^+))$$
and that by \eqref{26} we obtain $u\in C([0,T];H^3(\er^+))$.\\

\bigskip

We now sketch the proof of Lemma \ref{lema1}. First, we need to set the Cauchy problem (\ref{sistaltf}-\ref{idataf}) in the framework of real spaces, in order to apply a variant of Theorem 6 in \cite{K} (see also \cite{DFO} for a related result in the quasilinear local case). We introduce the new variables $F_1=Re(F)$, $F_2=Im(F)$, $u_1=Re(u)$, $u_2=Im(u)$ and, with $U=(F_1,F_2,v)$, $F_{1_0}=Re (F_0)$ and $F_{2_0}=Im (F_0)$, the Cauchy problem (\ref{sistaltf}-\ref{idataf}) can be written as follows:
\begin{equation}
\label{cplemma}
\left\{\begin{array}{llll}
U_t+A(U)U=g(t,U)\\
U(x,0)=(F_{1_0},F_{2_0},v_0)\in (H_0^1(\er^+))^2\times (H^2(\er^+)\cap H_0^1(\er^+)),
\end{array}\right.
\end{equation}
where 
$$A(U)=\left(\begin{array}{ccc}
0&\Delta&0\\
-\Delta&0&0\\
0&0&\Big(\int v^2dx\Big)\frac{\partial}{\partial x}
\end{array}\right)$$
and
$$g(t,U)=\left(\begin{array}{ccc}
2|u|^2F_2-(u_1^2-u_2^2)F_2+2u_1u_2F_1+bF_2v-b^2u_2|\tilde{u}|_x^2-abu_2\Big(\int v^2dx\Big)\\
2|u|^2F_1-(u_1^2-u_2^2)F_1+2u_1u_2F_2-bF_1v+b^2u_1|\tilde{u}|_x^2+abu_1\Big(\int v^2dx\Big)\\
-b|\tilde{u}|_x^2.
\end{array}\right)$$
Note that $A(U)$ and $g(t,U)$ are non-local.\\
We now set $X=(H^{-1}(\er^+))^2\times L^2(\er^+)$, $Y=(H_0^{1}(\er^+))^2\times (H^{2}(\er^+)\cap H_0^{1}(\er^+))$ and introduce the isomorphism 
$$S=I-\Delta\,:\,Y\to X.$$
Moreover, 
$$A\,:\,U=(F_1,F_2,v)\in W\to G(X,1,\beta),$$
where $W$ is an open ball of $Y$ centered at the origin with radius $R>0$ and $G(X,1,\beta)$ denotes the set of linear operators $D$ in $X$ such that $-D$ generates a $C_0$-semigroup $\{e^{-tD}\}$, with 
$$\|e^{-tD}\|\leq e^{\beta t},\,t\in[0,+\infty[,\,\beta=\frac 12 \sup_{x\in\er^+}\Big|av_x\int v^2dx\Big|\leq aR^2C(R),$$
$C(\cdot)$ continuous.\\
It is easy to check that $g$ verifies, for fixed $T>0$, $\|g(t,U)\|_Y\leq \lambda$, $t\in[0,T]$, $U\in C([0,T]; W)$, and for $U\in W$, $SA(U)S^{-1}=A(U)$. Moreover (cf. \cite{DFO} for a similar result), we have, for each pair $(U,U^*)$, $U=(F_1,F_2,v)$ and $U^*=(F_1^*,F_2^*,v^*)$ in $C([0,T];W)$,  
\begin{equation}
\label{28}
\|g(t,U)-g(t,U^*)\|_{L^1(0,T;X)}\leq c(T) \sup_{0\leq t\leq T}\|U(t)-U^*(t)\|_X,
\end{equation}
where $c(\cdot)$ is a continuous increasing function with $c(0)=0$. Finally, it is easy to prove that 
$$\|A(U_1)-A(U_2)\|_{\mathcal{L}(Y,X)}\leq c_1\|U_1-U_2\|_X,$$
for $U_1,U_2\in W$, $c_1$ not depending on $t\in[0,T]$.\\
Now, Lemma \ref{lema1} is a consequence of Theorem 6 in \cite{K}, adapted to the non-local operator $A$ and the non-local right-hand-side $g(t,U)$. where the local condition (7.7) is replaced by \eqref{28}, which is sufficient for the proof.
\section{A blow-up result}
In what follows, for $f\in L^1(\er_+)$, we will denote $\displaystyle \int_{\er^+} fdx$ by $\displaystyle \int fdx$ or simply by $\displaystyle \int f$ when there is no ambiguity.\\
We begin this section by presenting some invariants of the system \eqref{sistemain}:
\begin{Propriedade} Let $(u,v)$ a local solution of the Cauchy-Problem (\ref{sistemain}-\ref{idata}-\ref{bc}) as in the Theorem \ref{TEL}. Then,
 the mass
 \begin{equation}
  \label{mass}
{P}(t)=\int_{\er^+}|u(x,t)|^2dx,
  \end{equation}
the energy
 \begin{multline}
  \label{energy}
  \mathcal{E}(t)=\frac 12\int_{\er^+}|u_x(x,t)|^2dx+\frac 14\int_{\er^+}|u(x,t)|^4dx\\
 \end{multline}
 $$+\frac b2\int_{\er^+}v|u|^2dx+\frac a8\Big(\int_{\er^+}v^2dx\Big)^2$$
 and the momentum 
  \begin{equation}
  \label{momentum}
  \mathcal{M}(t)=\int_{\er^+}v^2(x,t)dx-2Im\int_{\er^+}u(x,t)\overline{u}_x(x,t)dx
 \end{equation}
are conserved: $\forall t\in[0,T[,\,\mathcal{P}(t)=\mathcal{P}(0), \mathcal{E}(t)=\mathcal{E}(0)\,\textrm { and }\mathcal{M}(t)=\mathcal{M}(0)$.
\end{Propriedade}
\noindent
{\bf Proof}\\
\noindent
Multiplying (\ref{sistemain})-i by $\overline{u}$ and integrating the imaginary part gives $\displaystyle \frac{d}{dt}\int |u|^2=0$.\\
Also, multiplying (\ref{sistemain})-i by $\overline{u}_t$ and integrating the real part yields
\begin{equation}
 \label{energia1}
-\frac 12\frac{d}{dt}\int |u_x|^2=\frac 14\frac{d}{dt}\int |u|^4+\frac b2\int v|u|^2_t.
\end{equation}
Using (\ref{sistemain})-ii,
$$\int v_t|u|^2=-a\int v^2\int |u|^2v_x=a\int v^2\int v(|u|^2)_x=$$
$$=-\frac ab\int v^2\int vv_t=-\frac a{4b}\frac d{dt}\Big(\int v^2\Big)^2.$$
Combining this identity with (\ref{energia1}),
$$-\frac 12\frac{d}{dt}\int |u_x|^2=\frac 14\frac{d}{dt}\int |u|^4+\frac b2\frac d{dt}\int v|u|^2+\frac a{8}\frac d{dt}\Big(\int v^2\Big)^2,$$
that is, $\displaystyle \frac d{dt}\mathcal{E}(t)=0$.

\bigskip

\noindent
Finally, we notice, by integrating by parts, that 
\begin{equation}
\label{eq2}
\frac {d}{dt}\int u\overline{u}_x=\int u\overline{u}_{xt}+\int u_t\overline{u}_x=2\int u_t\overline{u}_x.
\end{equation}
Now, by the Duhamel's formula for the nonlinear Schr\"odinger equation, we have
$$u(t)=e^{i\Delta t}u_0-i\int_0^te^{i\Delta(t-s)}[|u|^2u+buv](s)ds.$$
Hence, if
\begin{equation}
\label{33}
u_0\in H^2_0(\er^+)\cap H^3(\er^+)
\end{equation}
we derive, taking the $x$-derivative in the formula, that $u\in C([0,T];H_0^2(\er_+))$. Under this hypothesis, 
since $u_t=iu_{xx}-i|u|^2u-ibuv$, we obtain by integrating by parts that
$$2Im\int u_t\overline{u}_x=2Re\int u_{xx}\overline{u}_{x}-2Re\int |u|^2u\overline{u}_x-2b Re\int vu\overline{u}_x=-2b\int v|u|^2_x.$$
By (\ref{sistemain})-ii,
$$2Im\int u_t\overline{u}_x=-2b\int v(|u|^2)_x=2\int vv_t=\frac{d}{dt}\int v^2,$$
and, in view of (\ref{eq2}), it follows that $\displaystyle \frac{d}{dt}\mathcal{M}(t)=0$.\hfill$\blacksquare$

\bigskip

\bigskip

\noindent
Now, assume that
\begin{equation}
\label{35}
xu_0, x^{\frac 12}v_0\in L^2(\er^+)
\end{equation}
Using the multiplier $\mu_{\epsilon}(x)=e^{-\epsilon x}$ and letting $\epsilon\to 0^+$ it is easy to justify (see \cite{intro2bis} and \cite{C} for similar arguments) the following formal computations: 
$$\frac 12 \frac{d}{dt} \int x^2|u|^2=-Im \int x^2\overline{u}\frac{\partial^2u}{\partial x^2}=2Im\int x\frac{\partial u}{\partial x}\overline{u}$$
and
$$\begin{array}{llll}
\displaystyle\frac 12 \frac{d^2}{dt^2} \int x^2|u|^2dx&=&\displaystyle 2\frac{d}{dt}Im\int x\frac{\partial u}{\partial x}\overline{u}dx\\
&=&\displaystyle2Im\int x\frac{\partial^2u}{\partial x\partial t}\overline{u}dx+2Im \int x\frac{\partial u}{\partial x}\frac{\partial\overline{u}}{\partial t}dx\\
&=&\displaystyle -2Im\int\frac{\partial u}{\partial t}\overline{u}dx-2Im\int x\frac{\partial u}{\partial t}\frac{\partial \overline{u}}{\partial x}dx+2Im\int x\frac{\partial u}{\partial x}\frac{\partial \overline{u}}{\partial t}dx\\
&=&\displaystyle -2Im\int\frac{\partial u}{\partial t}\overline{u}dx-4Im\int x\frac{\partial u}{\partial t}\frac{\partial \overline{u}}{\partial x}dx\\
&=&\displaystyle -2Re \int\frac{\partial^2 u}{\partial x^2}\overline{u}dx+2Re\int |u|^4dx+2bRe\int v|u|^2dx\\
&&\displaystyle -4Re\int x\frac{\partial \overline{u}}{\partial x}\Big(\frac{\partial^2 u}{\partial x^2}-|u|^2u-buv\Big)dx\\
&=&\displaystyle 2\int |u_x|^2dx+2\int |u|^4dx+2b\int|u|^2vdx-2\int x\frac{\partial}{\partial x}|u_x|^2dx \\
&&+\displaystyle \int x\frac{\partial }{\partial x}|u|^4dx+2b\int x\frac{\partial}{\partial x}|u|^2vdx\\
&=&\displaystyle  4\int |u_x|^2dx+\int |u|^4dx+2b\int|u|^2vdx\\
&&\displaystyle -2\int x\Big(v_t+a\Big(\int v^2dx\Big)v_x\Big)vdx\\
&=& \displaystyle  4\int |u_x|^2dx+\int |u|^4dx+2b\int|u|^2vdx\\
&&\displaystyle -\frac{d}{dt}\int xv^2dx+a\Big(\int v^2dx\Big)^2.
\end{array}$$
Hence, we have
\begin{Lema}
\label{Lema31}
Under the assumptions \eqref{33} and \eqref{35}, the functions
$$t\to \int x^2|u|^2dx\qquad\textrm{ and }\qquad t\to \int xv^2dx$$
are respectively in $C^2([0,T])$ and $C^1([0,T])$, and
\begin{multline}
\label{36}
\frac{d^2}{dt^2}\Big[\frac 12\int x^2|u|^2dx+\int_0^t\int xv^2dxd\tau \Big]\\=4\int |u_x|^2+\int |u|^4dx+2b\int|u|^2vdx+a\Big(\int v^2\Big)^2.
\end{multline}
\end{Lema}
Under the supplementary assumption that $xu_0\in L^2(\er^+)$ we can also deduce that the function 
$$t\to \int x|u|^2dx$$
is in $C^1([0,T])$ and that
\begin{equation}
\label{37}
\frac{d}{dt}\int x|u|^2dx=-2Im\int u\frac{\partial\overline{u}}{\partial x}dx.
\end{equation}
Hence, we have, by  \eqref{37} and by the conservation of the momentum the following result:
\begin{Lema}
\label{Lema32}
Under the assumption \eqref{33} and and with $xu_0\in L^2(\er^+)$ we have, for all $t\in [0,T[$,
\begin{equation}
\label{38}
\frac{d}{dt}\int x|u|^2dx=\mathcal{M}(0)-\int v^2dx.
\end{equation}
\end{Lema}
Since, by \eqref{36} and the conservation of the energy, we have
$$\frac{d^2}{dt^2}\Big[\frac 12\int x^2|u|^2dx+\int_0^t\int xv^2dxd\tau \Big]=8\mathcal{E}(0)-\int |u|^4dx-2b\int |u|^2vdx,$$
we obtain that
$$\frac{d^2}{dt^2}\Big[\frac 12\int x^2|u|^2dx+\int_0^t\int xv^2dxd\tau+b^2\int_0^tx|u|^2dxd\tau \Big]$$
$$=8\mathcal{E}(0)-\int |u|^4dx-2b\int |u|^2vdx-b^2\int v^2dx+b^2\mathcal{M}(0)\leq 8\mathcal{E}(0)+b^2\mathcal{M}(0).$$
Hence, under the assumptions \eqref{33} and \eqref{35}, the function
$$\phi\,:\,t\to \frac 12\int x^2|u|^2dx+\int_0^t\int xv^2dxd\tau+b^2\int_0^t\int x|u|^2dxd\tau$$
is in $C^2([0,T])$, $\phi\geq 0$ and verifies the inequality
\begin{equation}
\label{39}
\frac{d^2}{dt^2}\phi(t)\leq 8\mathcal{E}(0)+b^2\mathcal{M}(0).
\end{equation}
By a classical virial argument, we have proven the following result:
\begin{Teorema}
\label{Teorema33} Under the assumptions of Theorem \ref{TEL} and \eqref{33},\eqref{35}, if
$$8\mathcal{E}(0)+b^2\mathcal{M}(0)<0,$$
then there is no global in time solution to the Cauchy Problem (\ref{sistemain}-\ref{idata}-\ref{bc}).\\
\end{Teorema}

\section{Existence of a global weak solution}
Let us fix the time interval $[0,T]$, $T>0$. We will use the approach known as the vanishing viscosity method and so we start by studying the following regularized version of (\ref{sistemain}):
\begin{equation}
 \label{sistemaep}
 \left\{\begin{array}{llll}
iu^{\epsilon}_t+u_{xx}=|u^{\epsilon}|^2u^{\epsilon}+bu^{\epsilon}v^{\epsilon}&(i)
\\
\displaystyle v^{\epsilon}_t+a\Big(\int_{\er^+}{v^{\epsilon}}^2dx\Big)v^{\epsilon}_x=-b(|u^{\epsilon}|^2)_x+\epsilon v^{\epsilon}_{xx}&(ii)
        \end{array}\right.
\end{equation}
with the initial conditions
\begin{equation}
\label{idataep}
\left\{\begin{array}{llll} u^{\epsilon}(x,0)=u_0(x),\quad u_0\in H_0^2(\er^+) \textrm{ and } x^{\frac 12}u_0\in L^2(\er^+),\\
 v^{\epsilon}(x,0)=v_0(x),\quad v_0\in H^2(\er^+)\cap H_0^1(\er^+).
 \end{array}\right.
\end{equation}
We start by stating the corresponding mass and energy identities: if
$$(u^{\epsilon},v^{\epsilon})\in\Big(C([0,T];H_0^1(\er^+))\Big)^2$$
is a solution of the above Cauchy problem, let
 \begin{multline*}
  \mathcal{E^{\epsilon}}(t)=\frac 12\int_{\er^+}|u_x^{\epsilon}(x,t)|^2dx+\frac 14\int_{\er^+}|u^{\epsilon}(x,t)|^4dx\\
  +\frac b2\int_{\er^+}v^{\e}|u^{\e}|^2dx+\frac a8\Big(\int_{\er^+}(v^{\e})^2dx\Big)^2\\
 \end{multline*}
 and
 \begin{equation*}
  \mathcal{M^{\epsilon}}(t)=\int_{\er^+}{v^{\epsilon}}^2(x,t)dx-2Im\int_{\er^+}u^{\epsilon}(x,t)\overline{u^{\epsilon}}_x(x,t)dx
 \end{equation*}
\noindent
Then, for all $t\in[0,T]$,

\begin{equation}
\label{43}
\|u^{\epsilon}(t)\|_2=\|u_0\|_{L^2}
\end{equation}
and
\begin{equation}
\label{44}
\mathcal{E}^{\epsilon}(t)+\epsilon \frac b2\int_0^t\int v^{\epsilon}_x(|u^{\epsilon}|^2)_xdxd\tau+\epsilon\frac a2\int_0^t\Big(\int (v^{\epsilon})^2dx\Big)\Big(\int (v_x^{\epsilon})^2dx\Big)d\tau=\mathcal{E}(0).
\end{equation}

\medskip

\noindent
Moreover, if $(u^{\epsilon},v^{\epsilon})\in C([0,T];H_0^2(\er^+))\times (C([0,T];H_2(\er^+)\cap H_0^1(\er^+))$ it is also easy to derive that, for all $t\in[0,T]$, 
\begin{equation}
\label{45}
\mathcal{M}^{\epsilon}(t)+2\epsilon\int_0^t\int (v^{\epsilon}_x)^2dxd\tau=\mathcal{M}(0).
\end{equation}
We point out that, since $u_0\in H_0^2(\er^+)$, if we prove the existence and uniqueness of a solution 
$$(u^{\epsilon},v^{\epsilon})\in \Big(C([0,T];H^2(\er^+)\cap H_0^1(\er^+))\Big)^2$$
to (\ref{sistemaep}-\ref{idataep}) we can derive, by the Duhamel formula for the first equation, that 
$u_\epsilon\in C([0,T];H_0^2(\er^+))$. Also, since $x^{\frac 12}u_0\in L^2(\er^+)$, we can also show that $x^{\frac 12}u^{\epsilon}\in C^1([0,T];L^2(\er^+))$ and that
\begin{equation}
\label{46}
\frac d{dt}\int x|u^{\epsilon}|^2dx=-2Im\int u^{\epsilon}\overline{u^{\epsilon}}_xdx.
\end{equation} 
Next, assume that we have proved the following result:
\begin{Lema}
\label{Lema41}
There exists $T'\in]0,T]$ and a unique
$$(u^{\epsilon},v^{\epsilon})\in \Big(C([0,T'];H^2(\er^+)\cap H_0^1(\er^+))\Big)^2$$
solution to the Cauchy problem (\ref{sistemaep}-\ref{idataep}) in $[0,T']$.
\end{Lema}
To prove that $T'=T$ we need the following Proposition:
\begin{Propriedade}
\label{Prop42} There exists $\epsilon_0>0$ and a continuous function $h$ (both independent of $a$) such that, for all $\epsilon\leq \epsilon_0$, if $(u^{\e},v^{\e})$ is a solution of (\ref{sistemaep}-\ref{idataep}), then, for all $t\in[0,T]$,
\begin{equation}
\label{47}
\|u^{\e}_x(t)\|_2\leq h(t)
\end{equation}
\end{Propriedade}
\noindent
{\bf Proof} In what follows, we drop the superscript $\e$.\\ Let $g(t)=\|u_x(t)\|_{L^2}$. Given $\delta, \delta_1>0$, by \eqref{43}, \eqref{44}, \eqref{45} and the Gagliardo-Nirenberg inequality, we deduce that
\begin{displaymath}
\begin{array}{llll}
\Big|\int v|u|^2dx\Big|&\leq& |b|\|v\|_2\|u\|_4^2\leq \delta\|v\||_2^4+c(\delta)\|u\|_4^{\frac 83}\\
&\leq&\delta\|v\||_2^4+c(\delta)\|u_0\|_2^2\|u_x\|_2^{\frac 23}\\
&\leq&\delta\|v\||_2^4+c(\delta)(\delta_1\|u_x\|_2^2+c(\delta_1)\|u_0\|_2^3)
\end{array}
\end{displaymath}
and 
$$\|v\|_2^4\leq c_0+c_1\|u_0\|_2^2\|u_x\|_2^2.$$
By choosing convenient $\delta$ and $delta_1$ we derive, by \eqref{44}, 
\begin{displaymath}
\begin{array}{llll}
g(t)=\|u_x(t)\|_2^2&\leq&\displaystyle c_2+c_2\e\int _0^t\int v_x(|u|^2)_xdxd\tau\\
&\leq&\displaystyle c_2+c_3\e\int_0^t \|v_x\|_2g(\tau)^{\frac 34}d\tau,
\end{array}
\end{displaymath}
with $c_2$ and $c_3$ independent of $\epsilon$ and $a$. Now, by \eqref{45}, we can proceed as in the proof of Lemma 1.3 in \cite{DF} to conclude the proof of Proposition 4.2.\hfill$\blacksquare$

\bigskip

Combining this result with \eqref{45} we derive, for $\epsilon\leq \epsilon_0$ (independent of $a$), 
\begin{equation}
\label{48}
\int (v^{\e})^2dx+\e\int_0^t(v_x^{\e})^2dxd\tau\leq h_1(t),\quad t\in[0,T],
\end{equation}
where $h_1$ is a continuous function independent of $\e$ and $a$.

\bigskip

\noindent
We pass to the proof that $T'=T$. We have, for $t\leq T'$ and by once more dropping the $\e$ superscript,
\begin{equation}
\label{49}
u(t)=e^{it\Delta}u_0-i\int_0^te^{i\Delta(t-\tau)}[|u|^2u+buv](\tau)d\tau
\end{equation}
and
\begin{equation}
\label{410}
v(t)=e^{t\Delta}v_0+\int_0^te^{i\Delta(t-\tau)}[-b(|u|^2)_x-a\Big(\int v^2dx\Big)v_x](\tau)d\tau
\end{equation}
By \eqref{43}, \eqref{47}, \eqref{48} and the well-known estimates for the heat equation
\begin{equation*}
\Big\|\frac{\partial}{\partial x}e^{t\Delta}\rho\Big\|_{\infty}\leq \frac c{t^{\frac 34}}\|\rho\|_1,\qquad \Big\|\frac{\partial}{\partial x}e^{t\Delta}\rho\Big\|_2\leq \frac c{t^{\frac 12}}\|\rho\|_2,\quad t>0,
\end{equation*}
we derive from \eqref{410} by the Gronwall's inequality that for $t\leq T'$
\begin{equation}
\label{411}
\|v_x(t)\|_2\leq h_2(t),
\end{equation}
where $h_2$ is continuous in $[0,T]$.\\
Now, from (\ref{sistemaep}-i), and since we have an estimate for $\Big\|b|u|^2+a\Big(\int v^2dx\Big)\Big\|_{H^1}$, we derive for $\|v_{xx}(t)\|_2$ an estimate similar to \eqref{411} using the properties of the heat kernel.\\
Finally, by \eqref{49}, the previous estimates and the properties of the Schr\"odinger group, we can deduce a similar estimate for $\|u_{xx}(t)\|_2$, and this achieves the proof that $T'=T$.

\bigskip

\bigskip
\noindent
We now sketch the prooof of Lemma \ref{Lema41}. To this end we consider the product space $\be_R^T\times B_R^T$, where $\be_R^T$ denotes the set of complex valued functions $\omega$ such that
$$\omega\in X_T=C([0,T];H^2(\er^+)\cap H^1_0(\er^+))\,:\,\|\omega\|_{L^{\infty}(0,T,H^2)}\leq R$$
and $B_R^T$ is the correspondent real-valued version of $\be_R^T$.\\
We consider the map 
$$(\tilde{u},\tilde{v})\in \be_R^T\times B_R^T\to (u,v)=\Phi(\tilde{u},\tilde{v})\in X_R^T\times X_R^T$$
where $(u,v)$ is the solution of the linear problem 
\begin{displaymath}
 \left\{\begin{array}{llll}
iu_t+u_{xx}=|\tilde{u}|^2\tilde{u}+b\tilde{u}\tilde{v}&
\\
\displaystyle v_t+a\Big(\int_{\er^+}{\tilde{v}}^2dx\Big)\tilde{v}_x=-b(|\tilde{u}|^2)_x+\epsilon \tilde{v}_{xx}
        \end{array}\right.
\end{displaymath}
with initial data (\ref{idataep}).\\
Since $$\displaystyle a\Big(\int_{\er^+}{\tilde{v}}^2dx\Big)\tilde{v}+b|\tilde{u}|^2\in C([0,T];H^2(\er^+)\cap H^1_0(\er^+))\cap C^1([0,T];L^2(\er^+))$$ it is not difficult to prove that there exists $T'\in]0,T]$ and $R>0$ such that the map $\Phi$ is a strict contraction of $\be_R^{T'}\times B_R^{T'}$ and, so, by the Banach fixed-point Theorem, we obtain a unique local in time solution of the Cauchy-Problem (\ref{sistemaep}-\ref{idataep}) in that space. Since one can take $T'=T$, as explained, and that, as observed above, if $u_0\in H_0^2(\er^+)$, $u\in C([0,T]; H_0^2(\er^+))$.\\ Hence, we have shown the following result:
\begin{Teorema}
\label{Teorema43}
There exists $\epsilon_0>0$ such that for each $T>0$ and $\e<\e_0$ there exists a unique solution 
$$(u^{\e},v^{\e})\in  C([0,T];H_0^2(\er^+))\times (C([0,T];H^2(\er^+)\cap H_0^1(\er^+))$$ of the Cauchy Problem (\ref{sistemaep}-\ref{idataep}). Moreover, $x^{\frac 12}u\in C([0,T];L^2(\er^+)).$
\end{Teorema} 
Now, for fixed $T>0$ and $\e<\e_0$, let $(u^{\e},v^{\e})$ be the solution obtained in Theorem \ref{Teorema43}. By Proposition \ref{Prop42} and \eqref{48}, we obtain that
\begin{equation}
\label{412}
\left\{\begin{array}{lllll}
\{u^\e\} \textrm{ is bounded in } L^{\infty}(]0,T[; H_0^1(\er^+))\\
\\
\displaystyle \Big\{ \frac{\partial u^\e}{\partial t}\Big\} \textrm{ is bounded in } L^{\infty}(]0,T[; H^{-1}(\er^+)),
\end{array}\right.
\end{equation}
\begin{equation}
\label{413}
\left\{\begin{array}{lllll}
\{v^\e\} \textrm{ is bounded in } L^{\infty}(]0,T[; L^2(\er^+))\\
\\
\displaystyle \Big\{ \frac{\partial u^\e}{\partial t}\Big\} \textrm{ is bounded in } L^{\infty}(]0,T[; H^{-2}(\er^+)),
\end{array}\right.
\end{equation}
and
\begin{equation}
\label{414}
\e\int_0^t\int (v^{\e}_x)^2dxd\tau\leq h_1(t),\,t\in[0,T],
\end{equation}
where $h_1$ is a continuous function independent of $\e$ or $a$. Recall that for each compact interval $I\subset\er^+$, the injection $H^1(I)\hookrightarrow L^2(I)$ is compact. Hence, the injection 
\begin{multline*}
W_I=\Big\{\omega\in L^2(]0,T[;H^1(\er^+))\,:\,\frac{\partial\omega}{\partial t}\in L^2(]0,T[;H^{-1}(\er^+))\Big\}\\
\hookrightarrow L^2(]0,T[, L^2(I))
\end{multline*}
is compact. By applying a standard diagonal extraction, there exists 
$$(u,v)\in L^{\infty}(]0,T[;H^1_0(\er^+))\times L^{\infty}(]0,T[;L^2(\er^+))$$ and a subsequence still denoted by $(u^\e,v^\e)_{\epsilon}$ such that
\begin{equation}
\label{415}
u^\e\rightharpoonup u \left\{\begin{array}{lll}
\textrm{ in } L^{\infty}(]0,T[; H_0^1(\er^+))\textrm{ weak}*,\\
\\
\textrm{ in } L^{2}(]0,T[; L^2(I))\textrm{ strong},\,I\subset\subset \er^+,\\
\\
\textrm{ a.e. in }]0,T[\times \er^+,
\end{array}\right.
\end{equation}
\begin{equation}
\label{415bis}
\frac{\partial u^\e}{\partial t}\rightharpoonup \frac{\partial u}{\partial t} 
\textrm{ in } L^{\infty}(]0,T[; H^{-1}(\er^+))\textrm{ weak}*,\\
\end{equation}
\begin{equation}
\label{416}
v^\e\rightharpoonup v 
\textrm{ in } L^{\infty}(]0,T[; L^2(\er^+))\textrm{ weak}*,\\
\end{equation}
\begin{equation}
\label{416bis}
\frac{\partial v^\e}{\partial t}\rightharpoonup \frac{\partial v}{\partial t} 
\textrm{ in } L^{\infty}(]0,T[; H^{-2}(\er^+))\textrm{ weak}*,\\
\end{equation}
hence
\begin{equation}
\label{417}
(u,v)\in C([0,T];L^2(\er^+))\times C([0,T];H^{-2}(\er^+)), \, u(0)=u_0\textrm{ and }v(0)=v_0. 
\end{equation}
Moreover, we have
$$\frac{d}{dt}\int (v^{\e})^2dx=2\int v^{\e}v_t^{\e}dx=-2b\int v^{\e}|u^{\e}|_x^2dx-2\e\int (v^{\e})^2dx,$$
and so, by \eqref{412}, \eqref{413} and \eqref{414}, for $e\leq \e_0$,
$$\int_0^t\Big|\frac{d}{dt}\int (v^{\e})^2dx\Big|dt\leq h_2(t),$$
where $h_2$ is a continuous function independent of $\e$ and $a$. Hence, $\Big\{\int (v^{\e})^2dx\Big\}_{\e}$ is bounded in $W^{1,1}(]0,T[)$ and, since the injection
$$W^{1,1}(]0,T[)\hookrightarrow L^q(]0,T[)$$
is compact for $1\leq q<+\infty$, we derive that, up to a subsequence,
\begin{equation}
\label{418}
\int (v^{\e})^2dx\to c(t)\,\,\textrm{ in }L^q(]0,T[)\,\textrm{ and a.e. in} ]0,T[.
\end{equation}
Note also that $c(t)\in L^{\infty}(]0,T[)$ (since $W^{1,1}(]0,T[)\hookrightarrow L^{\infty}(]0,T[)$). In particular, we deduce that 
\begin{equation}
\label{419}
\int (v^{\e}(x,t))^2dx\leq c(t)\,\textrm{ a.e. in }]0,T[.
\end{equation}
Recall that we have, by \eqref{45} and \eqref{46}, that
$$M(0)=\int v_0^2-2Im\int u_o\overline{u_0}_x=\int(v^{\e})^2-2 Im \int u^{\e}\overline{u^{\e}}_x+2\e\int_0^t\int (v^{\e}_x)^2dxd\tau$$
$$\geq \int (v^{\e})^2+\frac{d}{d t}\int x|u^{\e}|^2,$$
and so 
$$M(0)t+\int x|u_0|^2\geq \int_0^t c(\tau)d\tau+\liminf \int x|u^{\e}|^2.$$
It is easy to prove, by \eqref{46}, that $\displaystyle \int x|u|^2\leq \liminf \int x|u^{\e}|^2.$. Hence, we deduce, for $t\in [0,T]$, that
\begin{equation}
\label{420}
M(0)t+\int x|u|^2\geq \int_0^tc(\tau)d\tau+\int x|u|^2.
\end{equation}

By the previous estimates, we derive, letting $\e\to 0$ in \eqref{sistemaep} after a suitable integration,
\begin{multline}
 \label{421p17}
i\int_0^T\int u\phi_tdxdt+\int_0^T\int u_x\phi_xdxdt+\int u_0(x)\phi(x,0)dx\\
+\int_0^T\int |u|^2u\phi dxdt+b\int_0^Tvu\phi dxdt=0,
\end{multline}
\begin{multline}
 \label{422p17}
i\int_0^T\int v\psi_tdxdt+a\int_0^T c(t)\Big(\int v\psi_xdx\Big)dt+\int v_0(x)\psi(x,0)dx\\
-b\int_0^T\int (|u|^2)_x\psi dxdt=0,
\end{multline}
for all functions $\phi,\psi\in C_0^1(\er^+\times [0,T[)$, $\phi$ being complex and $\psi$ real-valued, and with
\begin{displaymath}
\begin{array}{llll}
u\in L^{\infty}(]0,T[;H_0^1(\er^+))\cap C([0,T];L^2(\er^+)),\\
\\
\displaystyle \frac{\partial u}{\partial t}\in L^{\infty}(]0,T[;H^{-1}(\er^+)),\\
\\
v\in L^{\infty}(]0,T[;L^2(\er^+))\cap C([0,T];H^{-1}(\er^+)),\\
\\
\displaystyle \frac{\partial v}{\partial t}\in L^{\infty}(]0,T[;H^{-1}(\er^+))
\end{array}
\end{displaymath}
since, in $\mathcal{D}'$,
\begin{equation}
 \label{423p17}
 \left\{\begin{array}{llll}
iu_t+u_{xx}=|u|^2u+buv
\\
\displaystyle v_t+ac(t)v_x=-b(|u|^2)_x
        \end{array}\right.
\end{equation}
and
\begin{equation}
 \label{425p17}
u(0)=u_0,\quad v(0)=v_0.
\end{equation}

Now, for each $\delta>0$ we choose a function $\gamma_{\delta}\in C^{\infty}(\er)$ such that
\begin{equation}
\label{426p17}
\gamma_{\delta}(x)=\left\{\begin{array}{lllll}
1\textrm{ for }x\geq \delta\\
0\textrm{ for }x\leq \frac{\delta}2,\quad \gamma_{\delta}\geq 0,\,\gamma_{\delta}'\geq 0.
\end{array}\right.
\end{equation}
We have the following auxiliary system in $\er^+\times]0,T[$ in $u$ and $v$ for each $\delta>0$:
\begin{equation*}
 \left\{\begin{array}{llll}
i(\gamma_{\delta}u)_t+\gamma_{\delta}u_{xx}=\gamma_{\delta}|u|^2u+b\gamma_{\delta}uv
\\
\displaystyle (\gamma_{\delta}v)_t+ac(t)\gamma_{\delta}v_x=-b\gamma_{\delta}(|u|^2)_x
        \end{array}\right.
\end{equation*}
that is
\begin{equation}
 \label{424}
 \left\{\begin{array}{llll}
i(\gamma_{\delta}u)_t+(\gamma_{\delta}u)_{xx}=2\gamma_{\delta}'u_x+\gamma_{\delta}''u+\gamma_{\delta}|u|^2u+b\gamma_{\delta}uv
\\
\displaystyle (\gamma_{\delta}v)_t+ac(t)(\gamma_{\delta}v)_x=ac(t)\gamma_{\delta}'v-b\gamma_{\delta}(|u|^2)_x.
        \end{array}\right.
\end{equation}
We can consider this system in $\er\times]0,T[$ by extending $u$ and $v$ to this set by $0$, with initial data ($u_0$ and $v_0$ also extended)
\begin{equation}
\label{425}
(\gamma_{\delta}u)(0)=\gamma_{\delta}u_0,\quad (\gamma_{\delta}v)(0)=\gamma_{\delta}v_0.
\end{equation}
We now consider, for $0<\e<\frac{\delta}2$, the usual family of mollifiers $\rho_{\e}\in \mathcal{D}(\er)$, with $0\leq \rho_{\e}\leq 1$, $supp\,\rho_{\e}\subset \overline{B}(0,\e)$ and $\displaystyle \int\rho_{e}=1$. We derive by \eqref{424} and \eqref{425} that
\begin{equation}
\label{426}
\left\{\begin{array}{lllll}
i(\rho_{\e}*\gamma_{\delta}u)_t+(\rho_{\e}*\gamma_{\delta}u)_{xx}=&2\rho_{\e}*\gamma_{\delta}'u_x+\rho_{\e}*\gamma_{\delta}''u+\rho_{\e}*\gamma_{\delta}|u|^2u\\
&+b\rho_{\e}*\gamma_{\delta}uv\\
(\rho_{\e}*\gamma_{\delta}v)_t+ac(t)(\rho_{\e}*\gamma_{\delta}v)_x=&ac(t)\rho_{\e}*\gamma'_{\delta}v-b\rho_{\e}*\gamma_{\delta}(|u|^2)_x,
\end{array}\right.
\end{equation}
with 
\begin{equation}
\label{427}
(\rho_{\e}*\gamma_{\delta}u)(0)=\rho_{\e}*\gamma_{\delta}u_0,\quad (\rho_{\e}*\gamma_{\delta}v)(0)=\rho_{\e}*\gamma_{\delta}v_0.
\end{equation}
Since $\frac{\delta}2-\e>0$, $\rho_{\e}*\gamma_{\delta}u$ and $\rho_{\e}*\gamma_{\delta}v$ are zero for $x\leq \frac{\delta}2-\e$.\\
Now, for simplicity, we drop the indexes $\e$ and $\gamma$ in the next computations. We easily derive that
\begin{equation}
\label{428}
\frac 12\frac{d}{dt}\int (\rho*\gamma v)^2=ac(t)\int (\rho*\gamma'v)(\rho*\gamma v)-b\int \rho*\gamma (|u|^2)_x(\rho*\gamma v).
\end{equation}
Passing to the limit in \eqref{428} when $\e\to 0$, we derive
\begin{equation}
\label{429}
\frac 12\frac{d}{dt}\int \gamma^2v^2=\frac a2c(t)\int (\gamma^2)'v^2-b\int\gamma^2(|u|^2_x)v.
\end{equation}
We also derive, from \eqref{428}, that
$$Im \frac d{dt}\int (\rho*\gamma u)(\rho*\gamma \overline{u})_x$$
and 
\begin{equation}
\label{430}
\begin{array}{lllll}
Im \int \frac{\partial}{\partial t}(\rho*\gamma u)(\rho*\gamma \overline{u})_x=&-&2Re\int (\rho*\gamma'u_x)(\rho*\gamma\overline{u})_x\\
&-&Re\int (\rho*\gamma''u)(\rho*\gamma\overline{u})_x\\
&-&Re\int (\rho*\gamma|u|^2u)(\rho*\gamma\overline{u})_x\\
&-&bRe\int (\rho*\gamma vu)(\rho*\gamma\overline{u})_x.
\end{array}
\end{equation}
We have, letting $\e\to 0$,
\begin{multline}
\label{431}
-2Re\int (\rho*\gamma'u_x)(\rho*\gamma\overline{u})_x\to -2Re\int \gamma'u_x(\gamma'\overline{u}+\gamma\overline{u}_x)\\
=-\int (\gamma^2)'|u_x|^2-\int(\gamma')^2(|u|^2)_x.
\end{multline}
In addition, we derive\\
$\displaystyle -Re\int (\rho*\gamma''u)(\rho*\gamma\overline{u})_x\to-Re\int \gamma''\gamma'|u|^2-Re\int \gamma''\gamma u\overline{u}_x$
\begin{equation}
 \label{432}
 \begin{array}{lllll}
&=&\displaystyle -\frac 12\int((\gamma')^2)_x|u|^2-\frac 12\int \gamma''\gamma(|u|^2)_x\\
\\
&=&\displaystyle\frac 12\int(\gamma')^2(|u|^2)_x+\frac 12\int(\gamma')^2(|u|^2)_x+\frac 14\langle (|u|^2)_{xx},(\gamma^2)'\rangle_{H^{-1}\times H^{1}_0}\\
\\
&=&\displaystyle\frac 12\int(\gamma')^2(|u|^2)_x+\frac 14\int \gamma'''|u|^2.
\end{array}
\end{equation}
Hence, by \eqref{431} and \eqref{432} we deduce
\begin{equation}
\label{433}
\lim_{\e\to 0}\Big(-2Re\int (\rho*\gamma'u_x)(\rho*\gamma\overline{u})_x-Re\int (\rho*\gamma''u)(\rho*\gamma\overline{u})_x\Big)
\end{equation}
$$=-\frac 12\int (\gamma^2)'|u_x|^2+\frac 14\int (\gamma^2)'''|u|^2.$$
Moreover, we have 
\begin{multline}
\label{434}
-Re\int (\rho*\gamma|u|^2u)(\rho*\gamma\overline{u})_x\to -\frac 12\int(\gamma^2)'|u|^4-Re\int \gamma^2|u|^2u\overline{u}_x\\=-\frac 14\int (\gamma^2)'|u|^4.
\end{multline}
Finally,
\begin{multline}
\label{435}
-bRe \int(\rho*\gamma vu)(\rho*\gamma\overline{u})_x\to-bRe\int \gamma vu(\gamma\overline{u})_x=\\-\frac 12b\int\gamma^2v(|u|^2)_x-\frac 12b\int(\gamma^2)'v|u|^2.
\end{multline}
Hence, by the previous estimates, we obtain
$$\frac 12\frac{d}{dt}\int \gamma^2v^2=Im \frac{d}{dt}\int \gamma u(\gamma\overline{u})_x+\frac a2 c(t)\int (\gamma^2)'v^2+\int (\gamma^2)'|u_x|^2$$
$$+ \frac 12\int (\gamma^2)'|u|^4+b\int(\gamma^2)'v|u|^2-\frac12\int(\gamma^2)'''|u|^2.$$ 
Since 
$$\Big|b\int (\gamma')^2|u|^2vdx\Big|\leq \frac 12 \int (\gamma')^2|u|^4+\frac{b^2}2\int (\gamma^2)'v^2$$
we derive, integrating in time,
$$\Big(\int \gamma^2v^2dx-2Im\int\gamma^2u\overline{u}_xdx\Big)-\Big(\int \gamma^2v_0^2dx-2Im\int\gamma^2u_0\overline{u_0}_xdx\Big)\geq $$
$$-b^2\int_0^t\int(\gamma')^2v^2dxd\tau-\int_0^t\int (\gamma^2)'''|u|^2dxd\tau,$$
and so, integrating over $[0,t]$ and letting $\delta\to 0$,
\begin{equation}
\label{436}
\int_0^t\Big(\int v^2dx-2Im\int u\overline{u}_xdx\Big)d\tau-M(0)t
\end{equation}
$$\geq \liminf_{\delta\to 0} \Big(-\int_0^t\int_0^{\tau}\Big[b^2\int(\gamma_{\delta}^2)'v^2dx+(\gamma_{\delta}^2)'''|u|^2\Big]dsd\tau \Big).$$
Since $x^{\frac 12}u_0\in L^2(\er^+)$ it is easy to derive that $x^{\frac 12}u\in W^{1,1}(]0,T[;L^2(\er^+))$ and 
$$\frac d{dt}\int x|u|^2=-2Im\int u\overline{u}_x.$$
Hence, by \eqref{435}, we deduce
\begin{equation}
\label{437}
\int_0^t\int v^2dxd\tau+\int x|u|^2dx\geq M(0)t+\int x|u_0|^2+
\end{equation}
$$\liminf_{\delta\to 0} \Big(-\int_0^t\int_0^{\tau}\Big[b^2\int(\gamma_{\delta}^2)'v^2dx+(\gamma_{\delta}^2)'''|u|^2\Big]dsd\tau \Big).$$
By \eqref{419}, \eqref{420} and \eqref{437} we deduce that $\displaystyle \int_0^t c(\tau)d\tau\geq \int_0^t\int v^2dxd\tau$ and
\begin{multline}
\label{438}
\int_0^t\int v^2dxd\tau\geq \int_0^tc(\tau)d\tau\\+\liminf_{\delta\to 0} \Big(-\int_0^t\int_0^{\tau}\Big[b^2\int(\gamma_{\delta}^2)'v^2dx+(\gamma_{\delta}^2)'''|u|^2\Big]dsd\tau\Big)
\end{multline}
a.e. in $[0,T]$.\\
We can state the following resut:
\begin{Teorema}
\label{T44}
Let $T>0$, $u_0\in H_0^2(\er^+)$, $x^{\frac 12}u_0\in L^2(\er^+)$ and $v_0\in H^2(\er^+)\cap H_0^1(\er^+)$.\\
Then there exists
$$u\in L^{\infty}(]0,T[;H_0^1(\er^+))\cap C([0,T]; L^2(\er^+))$$
and
$$v\in L^{\infty}(]0,T[;L^2(\er^+))\cap C([0,T]; H^{-1}(\er^+))$$
verifying the Cauchy problem \eqref{421p17}-\eqref{422p17} in $[0,T]$ and 
such that $$x^{\frac 12}u\in W^{1,1}(]0,T[;L^2(\er^+)),$$ $$\displaystyle \frac{\partial u}{\partial t}\in L^{\infty}(]0,T[;H^{-1}(\er^+))$$ and $$\displaystyle \frac{\partial v}{\partial t}\in L^{\infty}(]0,T[;H^{-1}(\er^+)).$$
Furthermore, there exists $c(t)\in L^{\infty}(]0,T[)$ satisfying
$$c(t)\geq \int v^2(x,t)dx$$
and \eqref{438} a.e. in $]0,T[$ with $\gamma_{\delta}$ defined in \eqref{426p17}.
\end{Teorema}
\begin{Rem}
The second term on the right-hand-side of \eqref{438} is a sort of boundary layer term, null for smooth $u$ and $v$. In this case,
$$c(t)=\int v^2(x,t)dx , \textrm {a.e. in }[0,T].$$
\end{Rem}
\bigskip
\noindent
\section{Bound states}
\noindent
In this section, we look for the existence of traveling-wave solutions to \eqref{sistemain} of the form
\begin{equation}
 \label{bound}
\left\{\begin{array}{llllll}
        u(x,t)&=&e^{i(\omega t +kx)}r(x+st)\\
v(x,t)&=&w(x+st),
       \end{array}\right.
\end{equation}
where $s\geq 0$ and $r,w$ are real-valued functions defined in $\er_0^+$.\\
By taking $k=\frac s2$ we obtain the system
\begin{equation}
 \label{sistema}
 \left\{\begin{array}{lll}
         -(\omega+\frac{s^2}4)r+r''=r^3+bwr\\
         sw+aw\int w^2=-br^2+c,\quad c\in\er,
        \end{array}\right.
\end{equation}
that is
\begin{equation}
 \label{sistema2}
 \left\{\begin{array}{lll}
        \mu r+r''=\lambda r^3\\
         \displaystyle w=\frac{-b}{a\alpha+s}\Big(r^2-\frac cb\Big),\quad\quad\alpha=\int w^2,
        \end{array}\right.
\end{equation}
where
$$\mu=-\omega-\frac{s^2}4-\frac{bc}{s+a\alpha}\,\textrm{ and }\,\lambda=1-\frac{b^2}{s+a\alpha}.$$

\bigskip

\noindent
Fixing $\lambda$, $\mu$ and $\gamma$, the general theory of ODEs states the existence of a maximal interval $I=[0,\overline{x}[$, $\overline{x}\in \er^+\cup\{+\infty\}$, such that the Initial Value Problem
\begin{equation}
 \label{ivp}
\left\{\begin{array}{lllll}
 \mu r+r''=\lambda r^3\\
r(0)=0,\,r'(0)=\gamma
       \end{array}\right.
\end{equation}
admits a unique smooth solution in $C^{\infty}([0,\overline{x}[)$. Furthermore, an elementary computation yields the following property:
 
\begin{Propriedade}
\label{cons1}
Let $r$ a smooth solution of (\ref{ivp}). Then for all $0\leq x<\overline{x}$,
\begin{equation}
\label{cons} 
E(x):=(r')^2(x)+\mu r^2(x)-\frac{\lambda}2r^4(x)=\gamma^2.
\end{equation}
\end{Propriedade}
In particular, note that if $\displaystyle \limsup_{x\to \overline{x}}r(x)=+\infty$ (i.e. $c=0$), then $\displaystyle\limsup_{+\infty} (r')^2=\gamma^2$ which is only possible if $\gamma=0$. Then, $r\equiv 0$.\\
\\
Next, we prove the following result.
\begin{Propriedade}
\label{ode1}
Let $\lambda,\mu>0$ and $\displaystyle \gamma=\frac{\mu}{\sqrt{2\lambda}}$.\\
Then \eqref{ivp} admits a unique solution $r\in C^{\infty}([0,+\infty[).$
Furthermore, $r$ is increasing and   $$r_{\infty}:=\lim_{x\to+\infty}r(x)=\sqrt{\frac{\mu}{\lambda}}.$$
\end{Propriedade}
\bigskip
\noindent
{\bf Proof}
Let us now consider the polynomial $P(x)=\frac {\lambda}2x^2-\mu x+\gamma^2$.\\ Note that we have, for all $x$, 
$$(r')^2(x)=P(r^2(x))\textrm{ and }P(r_0)=0\Leftrightarrow r_0=\frac {\mu}{\lambda}.$$
We begin by showing that $r'(x)> 0$ for all $x$. Indeed, let $$x_0=\min\{x>0\,:\,r'(x)=0\}.$$ Then, $\displaystyle r(x_0)=\sqrt{\frac{\mu}{\lambda}}$.\\
This is absurd by the uniqueness of the solution of  $\mu r+r''=\lambda r^3$ such that $r(x_0)=\sqrt{\frac{\mu}{\lambda}}$ and $r'(x_0)=0$. Indeed, note that the solution $\tilde{r}(x)\equiv \sqrt{\frac{\mu}{\lambda}}$ also satisfies these conditions.\\
 Hence, for all $x$, $r'(x)>0$ and for all $x<\overline{x}$, $r(x)<\sqrt{\frac{\mu}{\lambda}}$. This is enough to prove that $\overline{x}=+\infty$ and that $\displaystyle\lim_{x\to+\infty}r(x)=r_{\infty}.$\hfill$\blacksquare$

\bigskip

\bigskip

\noindent
We can now state the main result of this section.
 \begin{Teorema}   \label{teoremabs}
 Let $a,b\neq 0$. There exists $s, k, \omega$ such that (\ref{sistemain}) admits a smooth solution $(u,v)$ of the form \eqref{bound} such that
 \begin{itemize}
\item $r(0)=0$, $r\in L^{\infty}(\er^+)\cap C^{\infty}(\er^+)$ increasing;
\item $v\in L^2(\er^+)\cap C^{\infty}(\er^+)$;
 \end{itemize}

 \end{Teorema}
{\bf Proof} Let $a,b\neq 0$ and $\mu^*,\lambda^*>0$. We consider the solution $r(x)$ given by Proposition \eqref{ode1}.\\
 we put
$$w=\frac{-b}{s^*}\Big(r^2-\frac cb\Big)=\frac{-b}{s^*}(r^2(x)-r_{\infty}^2),$$
where we have chosen $\displaystyle c=\frac{b\mu}{\lambda}$ and $s^*$ such that $\lambda^*=1-\frac{b^2}{s^*}$.\\
We begin by checking that $w\in L^2(\er^+)$.\\
Setting $f=r_{\infty}-r$, $f\geq 0$, $f$ decreasing, and $\lim_{+\infty} f=0$. Furthermore, from Proposition \ref{cons1}, since $\gamma=\frac{\mu}{2\sqrt{\lambda}}$, we obtain that
$$(f')^2=-\mu(r^2-r_{\infty}^2)+\frac{\lambda}{2}(r^4-r_{\infty}^4)=(r-r_{\infty})(r+r_{\infty})(-\mu+\frac{\lambda}{2}(r^2+r_{\infty}^2))$$
$$=\frac{\lambda}{2}(r-r_{\infty})^2(r+r_{\infty})^2\geq \frac{\lambda}{2}(r-r_{\infty})^2r_{\infty}^2\geq \frac{\mu}2f^2,$$
where we have used twice that $\mu=\lambda r_{\infty}^2.$\\
Since $f'<0$, we obtain that $f'<-\sqrt{\frac\mu{2}}f.$ We can now infer that for $x\geq 0$
$$f(x)<f(0)e^{-\sqrt{\frac\mu{2}}x},$$
that is, $r\rightarrow r_{\infty}$ at exponential rate, yet $w\in L^2(\er^+)$.\\
Now, taking $s$ such that $s^*=a\int v^2+s$ and $\omega$ such that $\mu^*=-\omega-\frac{s^2}4-\frac{bc}{s^*}$, $(r,w)$ satisfies Theorem \ref{teoremabs}.
\hfill$\blacksquare$

\bigskip

\bigskip

\noindent
{\bf Acknowledgments}\\
The authors are grateful to Alain Haraux and Benedito Costa Cabral for stimulating discussions and suggestions. J.P. Dias acknowledges support of FCT (Portugal) grant UID/MAT/04561/2013. F. Oliveira  acknowledges support of FCT (Portugal) grant UID/MAT/00297/2013.

\bigskip

\small
\noindent \textsc{Jo\~ao-Paulo Dias}\\
CMAF-CIO and DM-FCUL \\
\noindent Campo Grande, Edif\'icio C6, 1749-016 Lisboa (Portugal)\\
\verb"jpdias@fc.ul.pt"

\bigskip

\small
\noindent \textsc{Filipe Oliveira}\\
Mathematics Department, FCT-UNL\\
Centro de Matem\'atica e Aplica\c{c}\~oes, CMA-UNL\\
NOVA University of Lisbon\\
Caparica Campus, 2829-516 (Portugal)\\
\verb"fso@fct.unl.pt"\\
\end{document}